\documentclass[12pt]{article}

\usepackage[utf8]{inputenc}

\usepackage{mathtools}
\usepackage{amsfonts,amsmath,amssymb,amsthm}

\usepackage[margin=2cm]{geometry}

\newtheorem{theorem}{Theorem}[section]

\title{On small non-uniform hypergraphs without property B}
\author{Danila Cherkashin\footnote{Institute of Mathematics and Informatics, Bulgarian Academy of Sciences. The work is supported by the Ministry of Education and Science of Bulgaria, Scientific Programme ``Enhancing the Research Capacity in Mathematical Sciences (PIKOM)'', No. DO1-67/05.05.2022}}

\begin{document}

\maketitle

\begin{abstract}
For a given hypergraph $H = (V,E)$ consider the sum $q(H)$ of $2^{-|e|}$ over $e \in E$.
Consider the class of hypergraphs with the smallest edge of size $n$ and without a 2-colouring without monochromatic edges. 
Let $q(n)$ be the smallest value of $q(H)$ in this class.

We provide a survey of the known bounds on $q(n)$ and make some minor refinements.
\end{abstract}

\textbf{Keywords:} non-uniform hypergraphs, hypergraph colouring, property B.

\textbf{MSC2020 classification:} 05C15, 05C65.

\section{Introduction}

A hypergraph $H = (V, E)$ is a finite set of vertices $V$ and a set of edges $E$ where
each edge is a set of at least two vertices. A 2-colouring of $H$ is an assignment
of colour blue or red to each vertex in $H$. A 2-colouring is \textit{proper} if each edge in $H$
is not monochromatic. We say that $H$ is \textit{2-colourable} if it admits a proper 2-colouring. 
A hypergraph is said to be \textit{$n$-uniform} if all of its edges have cardinality $n$. So a graph is just a $2$-uniform hypergraph.

A famous Erd{\H o}s -- Hajnal problem to find the minimum number of edges $m(n)$ in an $n$-uniform hypergraph that is not 2-colourable.
The best known asymptotic bounds are
\[
c \sqrt{\frac{n}{\ln n}} \leq m(n) \leq (1+o(1)) \frac{e \ln 2}{4} n^2 2^n,
\]
for a positive constant $c$ and $\ln$ standing for a natural logarithm. The lower bound was proved by Radhakrishnan and Srinivasan~\cite{radhakrishnan1998improved} and then another proof was given by Cherkashin and Kozik~\cite{cherkashin2015note}; the upper bound is due to Erd{\H o}s~\cite{Erdos2} and stays without improvements from 1963. A survey~\cite{raigorodskii2020extremalENG} is devoted to this problem and related topics.

Now let us pass to a non-uniform case. For a given hypergraph $H = (V,E)$ define the quantity
\[
q(H) := \sum_{e \in E} 2^{-|e|}.
\]
Note that $q(H)$ is the expectation of red edges in a random red-blue colouring of $V$ in which vertices get red colour with probability 1/2 independently on each other, and thus $q(H)$ is twice smaller than the expectation of monochromatic edges in such a colouring.

Erd{\H o}s~\cite{Erdos} asked in 1963 whether the function $q(n)$ is unbounded, where $q(n)$ is the minimal value of $q(H)$ over non-2-colourable hypergraphs $H$ with the minimal size of edge $n$. Beck~\cite{beck1980remark} in 1978 proved that $q(n) \geq  c\log^*t$ for $\log^*$ being the iterated logarithm and some positive constant $c$.
The in 2008 L. Lu~\cite{lu2008problem} announced a proof of a bound $q(n) \geq c \frac{\ln n}{\ln \ln n}$ but it turned out
to work only for simple hypergraphs (a hypergraph is \textit{simple} if every pair of edges shares at most 1 vertex). 
Shabanov~\cite{shabanov2015around} improved the lower bound for the class of simple hypergraphs to $c\sqrt{n}$ (also he made some refinements for hypergraphs with girth bounded from below in~\cite{shabanov2014coloring}). 
The best known asymptotic bounds on are
\[
c \ln n  \leq q(n) \leq (1+o(1)) \frac{e \ln 2}{4} n^2,
\]
where the lower bound is prove by Duraj, Gutowski, and Kozik~\cite{duraj2018note} and the upper bound is a direct consequence of the Erd{\H o}s upper bound on $m(n)$ and a straightforward estimate $q(n) \leqslant m(n)\cdot 2^{-n}$.

The main contribution of this note is a 2-time better asymptotic upper bound on $q(n)$.
\begin{theorem}
Let $n$ be an integer. Then
\[
q(n) \leq (1+o(1)) \frac{e \ln 2}{8} n^2.
\]
\label{mainth}
\end{theorem}

The proof is based on (probabilistic) alteration method and a resulting random hypergraph is the union of an $n$-uniform hypergraph and an $[n^2/4]$-uniform hypergraph. Akhmejanova~\cite{akhmejanova2021biuniform} refined some lower bounds on $q(n)$ for hypergraphs with only two different edge cardinalities.
Radhakrishnan and Srinivasan~\cite{radhakrishnan2021property} refined some lower bounds on $q(n)$ for the class of hypergraphs which locally have edges of comparable size.

\paragraph{Structure of the paper.} Section 2 is devoted to the case of small $n$. Section 3 contains the proof of Theorem~\ref{mainth}.

\section{The case of small $n$}

The values of $m(n)$ are known only for $n \leq 4$. We have $m(2) = 3$ with the only example of a triangle graph; and $m(3) = 7$ with the only example of the Fano plane. The best known upper bounds for small $n$ are reached by explicit examples; the current situation is outlined in Aglave, Amarnath, Shannigrahi and Singh~\cite{aglave2016improved}. No example of $q(n) < 2^{-n}m(n)$ is known. In this section we present an example of a hypergraph $H$ with the smallest edge with 4 elements such that
\[
2^{-4}m(4) < q(H) < 2^{-4}(m(4)+1),
\]
and the structure of $H$ completely differs from known examples for $m(4)$.

Now focus on $n=4$. In this case Seymour~\cite{seymour1974note} and Toft~\cite{toft1975color} independently showed that $m(4) \leqslant 23$. They used the example of a hypergraph on 11 vertices with the following edges:
\[
\{1,2,9,10\}, \quad \{3,4,9,10\}, \quad \{5,6,9,10\}, \quad \{7,8,9,10\},
\]
\[
\{1,2,9,11\}, \quad \{3,4,9,11\}, \quad \{5,6,9,11\}, \quad \{7,8,9,11\},
\]
\[
\{1,2,10,11\}, \quad \{3,4,10,11\}, \quad \{5,6,10,11\}, \quad \{7,8,10,11\},
\]
\[
\{1,3,5,8\}, \quad \{1,3,6,7\}, \quad \{1,4,5,7\}, \quad \{1,4,6,7\}, \quad \{1,4,6,8\},
\]
\[
\{2,3,5,7\}, \quad \{2,3,6,7\}, \quad \{2,3,6,8\}, \quad \{2,4,5,7\}, \quad \{2,4,5,8\}, \quad \{2,4,6,8\}.
\]
{\"O}sterg{\aa}rd~\cite{ostergaard2014minimum} by a complicated computer search show that $m(4)=23$, and there is only one example on at most 11 vertices.

We provide an example of a hypergraph $H$ with sixteen vertices, twenty 4-edges and sixty 8-edges which is not 2-colourable. 
It means that $q(H) = \frac{95}{64} = \frac{23}{16} + \frac{3}{64} < \frac{24}{16}$ which means that this is better than any known $4$-graph, except the Seymour -- Toft graph.

\begin{figure}[h]
    \centering
    \includegraphics[width=0.6\textwidth]{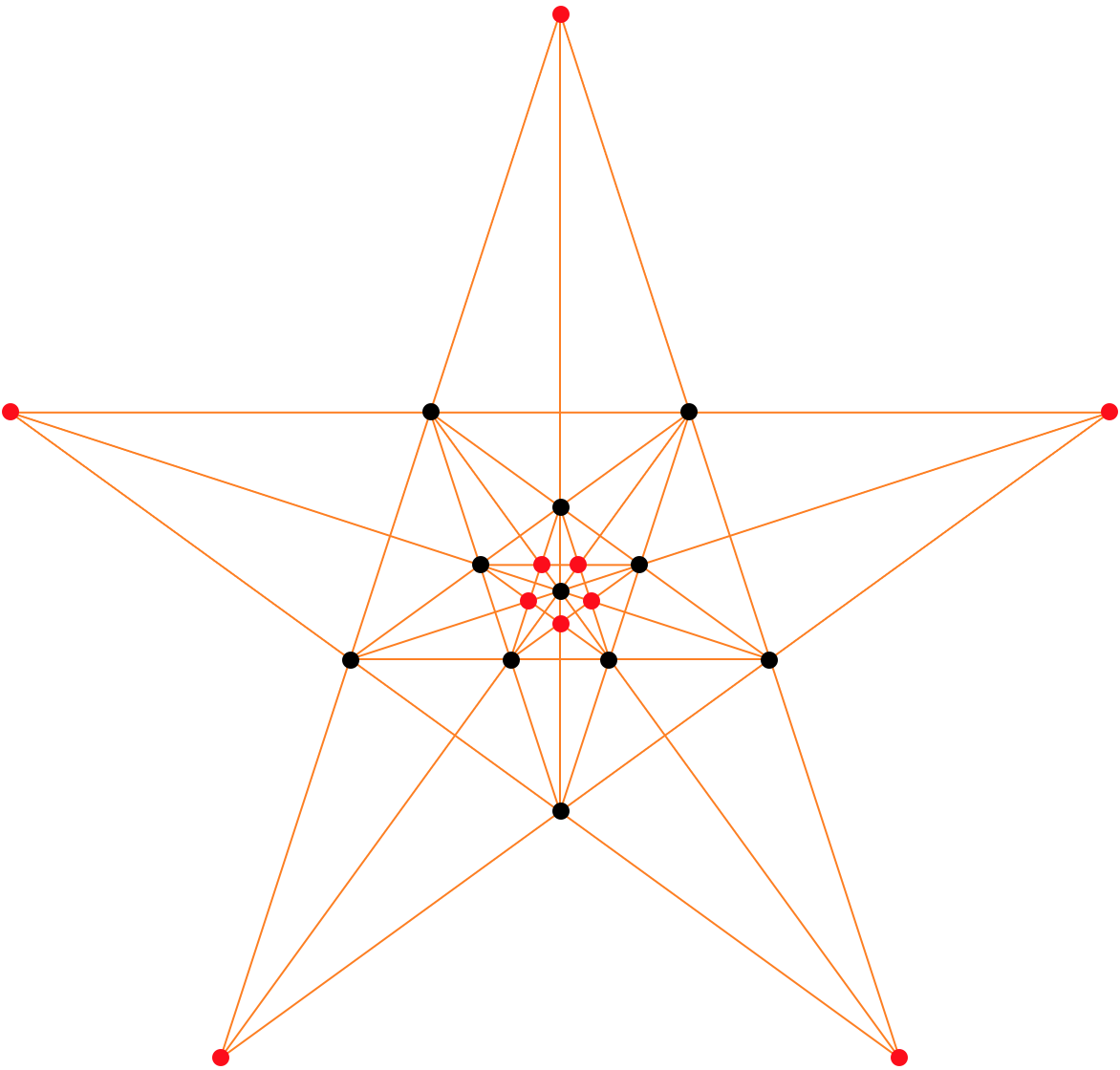}
    \caption{An explicit definition of $H_4$. Points are vertices and lines are edges. In this picture, opposite red points are the same. For example, the vertical line appears to contain 5 points, but the 2 red points are the same, so it just contains 4 points}
    \label{Affine4}
\end{figure}

By construction $H$ consists of a 4-uniform part $H_4$ and an 8-uniform part $H_8$.
The 4-uniform part $H_4$ is an affine plane over $GF(4)$. It may be also defined in an explicit way, see Fig.~\ref{Affine4} (I found out this representation from the answer of user Matt at math.stackexchange.com). A direct computation shows that $H_4$ has $120$ proper 2-colourings and every proper 2-colouring has exactly 8 red and 8 blue vertices.\footnote{K. Vorob'ev told me that blue sets of those colourings form a 3-(16,8,12) design.}
Thus these colourings form 60 ``opposite'' pairs, id est colourings in a pair get from each other by swapping the colours.
Then taking a red vertices of one member from each pair as an 8-edge one get $|E(H_8)| = 60$ and $H_4 \cup H_8$ has no proper 2-colouring as desired.

\section{Proof of Theorem~\ref{mainth}}

To avoid rounding in calculations, assume that $n$ is even; the case of odd $n$ is analogues.
Consider a set of vertices $V$ with the cardinality $v = n^2/2$, and choose $m$ random edges uniformly and independently; the number $m$ will be specified later. Fix a colouring $C$; clearly the probability of the event that a randomly chosen edge is monochromatic is equal to
\[
p := \frac{\binom{v_1}{n} + \binom{v_2}{n}}{\binom{v}{n}},
\]
where $v_1$ and $v_2$ denote the numbers of vertices of the first and second colour,
respectively. Hence, since the edges are chosen independently, the probability is
$(1-p)^m$ that after choosing $m$ random independent edges the colouring $C$ is proper.
Put
\[
q :=  \frac{2\binom{v/2}{n}}{\binom{v}{n}}.
\]
It is well-known that
\[
q = (1+o(1)) \frac{2}{e \cdot 2^n}.
\]

Note that $p \geq q$ because of the convexity of the sequence $\left\{\binom{t}{n} \right\}_{t \geq 0}$. Since the
total number of colourings is $2^{n^2/2}$, and the probability that a fixed colouring is proper is bounded by $(1-q)^m$, the expectation of the number of proper colouring is at most
\[
2^{n^2/2} (1-q)^m < e^{\ln 2 \cdot n^2/2 - qm},
\]
we use a standard inequality $1-t < e^{-t}$, $t > 0$.

To get the Erd{\H o}s upper bound one should take
\[
m = (1 + o(1)) \frac{e \ln 2}{4}n^2 2^n
\]
and check that such a choice leads to $\ln 2 \cdot n^2/2 - qm < 0$ which means that with a positive probability a random graph with $m$ edges has no proper 2-colouring.

For our purpose we need a twice smaller number of edges, i.e. $m' = m/2$; then the expectation of the number of proper colouring is at most $2^{n^2/4}$, so with a positive probability a random graph $H_1$ with $m'$ edges of size $n$ has at most $2^{n^2/4}$ proper colourings.
For each such a colouring $C$ we consider an edge $e_C$ of size $n^2/4$ which is monochromatic in $C$.
Let $H_2$ be a hypergraph on the same vertex set $V$ consisting of all such edges $e_C$, and let $H$ be the union of $H_1$ and $H_2$.
Then 
\[
q(H) = q(H_1) + q(H_2) \leq (1+o(1)) \frac{e \ln 2}{8}n^2 + 1 = (1+o(1)) \frac{e \ln 2}{8}n^2,
\]
as desired.

\paragraph{Acknowledgments.} I am grateful to Alexey Gordeev for an independent computer search and to Konstantin Vorob'ev for his remarks.

\bibliographystyle{plain}
\bibliography{main}

Danila Cherkashin

Institute of Mathematics and Informatics

Bulgarian Academy of Sciences

Acad. G. Bonchev Str., Bl. 8

1113 Sofia, Bulgaria

e-mail: jiocb@math.bas.bg

\end{document}